\documentclass[10pt]{article}
\usepackage{amsfonts,amsmath,amsthm,amssymb}

\usepackage[vcentering,dvips]{geometry}
\newtheorem{thm}{Theorem}

\theoremstyle{remark}

\theoremstyle{definition}

\newcommand{\<}{\langle}                                                                
\renewcommand{\>}{\rangle}                                                              

\newcommand{\R}{\mathbb{R}}

\newcommand{\N}{\mathbb{N}}

\title{Weak and quasi-polynomial tractability of approximation of infinitely differentiable functions}

\author{Jan Vyb\'\i ral\footnote{Technical University Berlin, Department of Mathematics, Secretary office MA 4-1, Street of 17. June 136, 10623 Berlin, Germany, e-mail: {\tt vybiral@math.tu-berlin.de}}
\thanks{Research supported by the DFG Research Center {\sc Matheon} ``Mathematics for key technologies'' in Berlin.}}

\begin{document}

\maketitle

\centerline{\emph{Dedicated to J. F. Traub and G. W. Wasilkowski}}
\centerline{\emph{on the occasion of their 80th and 60th birthdays}}

\vskip1cm

\begin{abstract}
We comment on recent results in the field of information based complexity, which state (in a number of different settings),
that approximation of infinitely differentiable functions is intractable and suffers from the curse
of dimensionality. We show that renorming the space of infinitely differentiable functions in a suitable
way allows weakly tractable uniform approximation by using only function values. Moreover, the approximating algorithm is based on
a simple application of Taylor's expansion about the center of the unit cube. We discuss also the approximation on the Euclidean ball
and the approximation in the $L_1$-norm.
\end{abstract}

\noindent{\bf Key words:} weak tractability, uniform approximation, infinitely differentiable functions, curse of dimensionality

\section{Introduction}

We consider different classes $F_d$ of infinitely-differentiable functions $f:\R^d\to\R$
and discuss algorithms using only function values of $f$ in order to approximate $f$
uniformly, or in the $L_1$-norm. We are especially interested in the case of large $d\gg 1$.

In the classical setting of approximation theory, the dimension of the Euclidean space $d$ is fixed. Furthermore,
the decay of the minimal error $e(n)$ of approximation of smooth functions
in a Lebesgue space norm is very well studied for both algorithms using $n$ arbitrary linear functionals and
for algorithms using only $n$ function evaluations. We refer to \cite{K,N,NT,P,T,Vyb2,Vyb1} and references therein.
The decay is usually polynomial, speeds up with increasing smoothness and slows down with increasing dimension.
Furthermore, this terminology typically hides the dependence
of the constants on the dimension $d$, which might be even exponential. This motivates the question,
what happens if both the dimension $d$ and the smoothness parameter $s$ tend to infinity.

If $n(\varepsilon,d)$ denotes the minimal number of function values needed to approximate
all functions from $F_d$ up to the error $\varepsilon>0$, we say, that the problem suffers from
\emph{curse of dimensionality}, if $n(\varepsilon,d)$ grows exponentially in $d$. This means, that
there are positive numbers $c, \varepsilon_0$ and $\gamma$, such that
$$
n(\varepsilon,d)\ge c (1+\gamma)^d\quad\text{for all}\quad0<\varepsilon\le\varepsilon_0\ \quad\text{and infinitely many}\quad d\in\N.
$$
On the other hand, we say that the problem is \emph{weakly tractable} if
$$
\lim_{\varepsilon^{-1}+d\to\infty}\frac{\ln n(\varepsilon,d)}{\varepsilon^{-1}+d}=0.
$$
Furthermore, the problem is \emph{quasi-polynomially tractable} in the sense of \cite{GW}, if there exist two constants $C,t>0$, such that
\begin{equation}\label{eq:QPT}
n(\varepsilon,d)\le C\exp\left\{t(1+\ln(1/\varepsilon))(1+\ln(d))\right\}
\end{equation}
for all $0<\varepsilon<1$ and all $d\in\N.$ For the sake of completeness, we add that a problem is \emph{polynomially tractable} if there exist non-negative numbers
$C, p$ and $q$ such that
$$
n(\varepsilon, d) \le C \varepsilon^{-p} d^q\quad\text{for all}\quad 0<\varepsilon<1\quad\text{and}\quad d \in \N.
$$
If $q=0$ above then a problem is \emph{strongly polynomially tractable}.
We refer to the monographs \cite{Traub, NW1} and \cite{NW2} for a detailed discussion of these and other kinds of \emph{(in)tractability}
and the closely related field of \emph{information based complexity}.

The $L_\infty$-approximation of infinitely differentiable functions was studied in \cite{HZ}, where the authors showed, that the problem is not strongly polynomially tractable.
It was also discussed in \cite{NW1}, cf. Open Problem 2 therein. An essential breakthrough was achieved in \cite{NW_2009_2} (which in turn is based on \cite{NW_2009_1} and answers an open problem posed there),
where uniform approximation of the functions from the class
$$
{\mathbb F}_d=\{f:[0,1]^d\to\R:\sup_{\alpha\in\N_0^d}\|D^\alpha f\|_\infty\le 1\}
$$
was shown to satisfy $n(\varepsilon,d)\ge 2^{\lfloor d/2\rfloor}$ for all $0<\varepsilon<1$ and all $d\in\N$ and this result is also true if
arbitrary linear functionals are allowed as the information map about $f$.
Hence, the problem is intractable and suffers from the curse of dimensionality. In the context of weighted spaces of infinitely differentiable functions,
the problem was also discussed in \cite{Weimar}.

Multivariate integration of infinitely differentiable functions from the class ${\mathbb F}_d$ was conjectured not to be polynomially tractable in \cite{W}
and was shown not to be strongly polynomially tractable in \cite{Woj}. Furthermore, it is known (cf. \cite{Suk} and \cite{HNUW}) that multivariate integration
of functions from
$$
C^k_d=\{f:[0,1]^d\to\R:\sup_{\alpha:|\alpha|\le k}\|D^\alpha f\|_\infty\le 1\}
$$
suffers from the curse of dimensionality for all $k\in\N.$ 
Although multivariate integration of infinitely differentiable functions is also discussed in \cite{HNUW} and \cite{HNUW2}, it seems to be 
still an open problem if the curse of dimensionality holds also for multivariate integration and the class ${\mathbb F}_d$.

The main result of this paper is the following.
\begin{thm}\label{thm1.1}
\begin{itemize}
\item[(i)] Uniform approximation on the cube $[-1/2,1/2]^d$ of functions from the class
\begin{equation}\label{eq:Fd1}
F_d^1=\Bigl\{f\in C^\infty([-1/2,1/2]^d):\sup_{k\in\N_0}\sum_{|\beta|=k} \frac{\|D^\beta f\|_\infty}{\beta!}\le 1\Bigr\}
\end{equation}
is quasi-polynomially tractable.
\item[(ii)] Uniform approximation on the balls $B(0,r_d)=\{x\in\R^d:\|x\|_2\le r_d\}$, where $r_d$ are chosen in such a way, that the volume of $B(0,r_d)$ is equal to one,
and the functions are from the class
\begin{equation}\label{eq:Fd2}
F_d^2=\{f\in C^\infty(B(0,r_d)):\sup_{k\in \N_0} \|\partial^{k}_\nu f\|_\infty\le 1\}
\end{equation}
is weakly tractable. Here, $(\partial_\nu^k f)(x)$ denotes the $k$-th derivative of $f$ at $x\not=0$ in the ``normal'' direction $x/\|x\|_2$.
\item[(iii)] Approximation in the $L_1$-norm on $B(0,r_d)$, where
$r_d$ are as above and the functions are from the class
\begin{equation}\label{eq:Fd3}
F_d^3=\Bigl\{f\in C^\infty(B(0,r_d)):\sup_{k\in\N_0} \int_0^{r_d} S(\partial_\nu^{k}f,r) dr\le 1\Bigr\}
\end{equation}
is weakly tractable. Here $S(\partial_\nu^{k}f,r)$ are the averages of $|\partial_\nu^{k}f|$ on the sphere $r{\mathbb S}^{d-1}$.
\end{itemize}
\end{thm}

Our method is rather simple and involves only the Taylor's expansion of a smooth function about the center
of the domain under consideration. The next two sections of this paper are devoted to the proof of this theorem.
Finally, the last section presents some possible extensions of this method and mentions several closely related open problems.

\section{Uniform approximation}\label{Sec2}

We study first the uniform approximation of infinitely differentiable functions $f$ using only its function values.
Let us recall, that the paper of Novak and Wo\'zniakowski \cite{NW_2009_2} shows intractability of this problem
on the unit cube even for arbitrary linear functionals as the information map of $f$.
We show that modifying the norm in a suitable way
leads immediately to weak (and even quasi-polynomial) tractability, even when allowing only function values of $f$ as the admissible information and using only Taylor's expansion about the center of the cube or ball, respectively.
Of course, any tractability result for approximation using function values implies the same result also for algorithms using arbitrary linear information.
On the other hand, we leave it as an open problem if algorithms using general linear information could achieve even better art of tractability.

Let us recall the standard multivariate notation, which we shall use with connection to the Taylor's theorem. 
If $\alpha=(\alpha_1,\dots,\alpha_d)\in\N_0^d$ is a multiindex, we denote
\begin{align*}
|\alpha|&=\alpha_1+\dots+\alpha_d,\\
D^\alpha f(x)&=\frac{\partial^{|\alpha|}f(x)}{\partial x_1^{\alpha_1}\dots\partial x_d^{\alpha_d}},\quad x\in\R^d,\\
\alpha!&=\alpha_1!\dots\alpha_d!,\\
x^\alpha&=x_1^{\alpha_1}\dots x_d^{\alpha_d},\quad x\in\R^d.
\end{align*}

\subsection{Unit cube}\label{Sec2.1}
In this part, we study uniform approximation of infinitely differentiable functions on $[-1/2,1/2]^d.$ Our algorithm is based on Taylor's formula
\begin{align*}
f(x)&=(T_kf)(x)+\sum_{|\beta|=k+1} (R_\beta f)(x)x^\beta,
\end{align*}
where
\begin{align*}
(T_k f)(x)&=\sum_{|\alpha|\le k}\frac{(D^\alpha f)(0)}{\alpha!}x^\alpha,\\
(R_\beta f)(x)&=\frac{|\beta|}{\beta!}\int_0^1 (1-t)^{|\beta|-1}(D^\beta f)(tx)dt.
\end{align*}
The approximation of $f$ by $T_k f$ at the point $x\in [-1/2,1/2]^d$ results into an error
\begin{align*} 
|f(x)-T_kf(x)|&\le \sum_{|\beta|=k+1} |R_\beta f(x)|\cdot |x^\beta|\\
&\le \sum_{|\beta|=k+1} \frac{(k+1)\cdot |x^\beta|}{\beta!}\int_0^1 (1-t)^{k}|(D^\beta f)(tx)|dt\\
&\le \left(\frac{1}{2}\right)^{k+1}\sum_{|\beta|=k+1} \frac{(k+1)}{\beta!}\int_0^1 (1-t)^{k}dt\cdot \|D^\beta f\|_\infty\\
&=\left(\frac{1}{2}\right)^{k+1}\sum_{|\beta|=k+1} \frac{\|D^\beta f\|_\infty}{\beta!}.
\end{align*}
Hence,
\begin{align*}
\sup_{f\in F_d^1}\|f-T_kf\|_\infty\le \left(\frac{1}{2}\right)^{k+1}
\end{align*}
for every $k\in\N_0$, where $F_d^1$ was defined in \eqref{eq:Fd1}.
To discuss the tractability of the problem, we need to estimate also the number of points needed to recover $T_kf.$
As we are allowed to take only samples of $f$, and not of its derivatives, we are actually not able to recover $T_kf$ exactly.
But using finite order differences, we may approximate it to an arbitrary precision using a bounded number of points.

To estimate the number of sampling points needed to approximate all derivatives up to the order $k$, we use induction. To evaluate $f(0)$ we need to sample 
$f$ at the point $A_0=\{0\}$. To calculate the first order differences, we need the points from the set
$$
A_1=A_0\cup \bigcup_{j=1}^d (A_0+he_j),
$$
which has $d+1$ points. For the second order differences, we need the values at 
$$
A_2:=A_1\cup \bigcup_{j=1}^d (A_1+he_j)=\{h(\alpha_1,\dots,\alpha_d):\alpha\in\N_0^d,|\alpha|\le 2\},
$$
which has $\binom{d+2}{2}$ points. By induction we obtain, that to evaluate all finite order differences up to the order $k$, we need the values at
$$
A_k=\{h(\alpha_1,\dots,\alpha_d):|\alpha|\le k\}
$$
and that this set has
\begin{equation*}
\sum_{j=0}^k\binom{d+j-1}{j}=\binom{d+k}{k}
\end{equation*}
points.

To show that the problem is quasi-polynomial tractable, we proceed in the following way. Given an $1>\varepsilon>0$, we choose the smallest $k\in\N_0$, such that
$$
\varepsilon\ge \Bigl(\frac{1}{2}\Bigr)^{k+1},
$$
i.e. 
$$
k+1:=\left\lceil \frac{\ln(1/\varepsilon)}{\ln 2} \right\rceil,
$$
where $\lceil a\rceil$ denotes the smallest integer, which is larger than or equal to the real number $a\in\R.$
Together with the estimate
$$
n(\varepsilon,d)\le \binom{d+k}{k}\le \Bigl(\frac{e(d+k)}{k}\Bigr)^k,
$$
this gives that there is $t>0$, such that
$$
\ln n(\varepsilon,d)\le k(1+\ln(d+k)-\ln k)\le t(1+\ln(1/\varepsilon))(1+\ln d).
$$
This immediately implies \eqref{eq:QPT}, i.e. the quasi-polynomial tractability.

\subsection{Euclidean ball}\label{sec:eucl}

In this section, we discuss the uniform approximation of an infinitely differentiable function $f$ on an Euclidean ball in $\R^n$ with radius $r_d>0$.
As all the sets of infinitely differentiable functions under consideration include all constant functions with values between -1 and 1, the initial error
of approximation is always 1 for every sequence $(r_d)_{d\in\N}.$ Nevertheless, we shall be interested at most in the case, when
the Lebesgue measure of $B(0,r_d):=\{x\in\R^d:\|x\|_2\le r_d\}$ is one. It is very well known, cf. \cite{HNUW}, that $r_d\approx \sqrt{d}$ in this case, i.e. there are two absolute constants
$C>c>0$, such that $c\sqrt{d}\le r_d\le C\sqrt{d}.$

We start again with Taylor's expansion of an infinitely differentiable function $f$ about zero. Let $x\in B(0,r_d), x\not=0$ and let $g_x(t)=f(tx), 0\le t\le 1$.
Then
\begin{equation}\label{eq:tay2}
f(x)=g_x(1)=\sum_{j=0}^k\frac{g_x^{(j)}(0)}{j!}+\frac{1}{k!}\int_0^1(1-t)^kg_x^{(k+1)}(t)dt.
\end{equation}
Iterating the formula
$$
g_x'(t)=\<(\nabla f)(tx), x\>=\<(\nabla f)(tx), \frac{x}{\|x\|_2}\>\|x\|_2=(\partial_\nu f)(tx)\cdot\|x\|_2
$$
we obtain
$$
g_x^{(k+1)}(t)=(\partial^{k+1}_\nu f)(tx)\cdot \|x\|_2^{k+1},\qquad 0<t<1,
$$
where $(\partial_\nu f)(x)$ denotes the derivative of $f$ at $x\not=0$ in the direction $x/\|x\|_2$.

We denote
$$
\widetilde T_k f(x)=\sum_{j=0}^k\frac{g_x^{(j)}(0)}{j!},
$$
which allows to estimate the error of approximation of $f(x)$ by the Taylor's polynomial about $x$ by
\begin{align*}
|f(x)-\widetilde T_kf(x)|&\le \frac{\|x\|_2^{k+1}}{k!}\int_0^1 (1-t)^k\cdot|(\partial^{k+1}_\nu f)(tx)|dt\\
&\le \frac{r_d^{k+1}\|\partial_\nu^{k+1} f\|_\infty}{(k+1)!}\le \frac{r_d^{k+1}}{(k+1)!},
\end{align*}
if $f\in F_d^2$, where $F_d^2$ is as in \eqref{eq:Fd2}.

Finally, we observe that we need again $\binom{d+k}{k}$ points to approximate the derivatives $D^\alpha f(0)$ for every $|\alpha|\le k$, which again allows to approximate
all the derivatives $g_x^{(j)}(0)$ for all $x\in B(0,r_d)$ and all $0\le j\le k.$

Hence, if $d\in\N$ is fixed and $1>\varepsilon>0$ is given, we choose first the smallest $k$, for which
\begin{equation}\label{eq:num2}
\varepsilon \ge\frac{r_d^{k+1}}{(k+1)!}.
\end{equation}
This is always possible, as the right hand side goes to zero for $d$ fixed and $k\to\infty$. On the other side, let us mention that if \eqref{eq:num2} holds for 
$\varepsilon<1$, then $k$ is at least of the order $r_d.$  Using the estimate
$$
\Bigl(\frac{e r_d}{k+1}\Bigr)^{k+1}\ge\frac{r_d^{k+1}}{(k+1)!}
$$
we obtain that \eqref{eq:num2} is satisfied any time we have
$$
k\ge \max(e^2r_d, \ln(1/\varepsilon)).
$$
If finally $r_d\approx \sqrt{d}$, this implies the weak tractability of the problem by
$$
\lim_{\varepsilon^{-1}+d\to\infty}\frac{\ln n(\varepsilon,d)}{\varepsilon^{-1}+d}
\le \lim_{\varepsilon^{-1}+d\to\infty}\frac{k(1+\ln(d+k)-\ln k)}{\varepsilon^{-1}+d}
\le \lim_{\varepsilon^{-1}+d\to\infty}\frac{k\ln(d+k)}{\varepsilon^{-1}+d}=0.
$$
Unfortunately, the calculation above does not give quasi-polynomial tractability in this case.

\section{Approximation in the $L_1$-norm}\label{sec:numer}

In this section we prove the last part of Theorem \ref{thm1.1}.
Of course, on a domain with volume one, the error of approximation in the $L_1$-norm may be bounded from above by the
error of uniform approximation. Therefore, the problem is weakly tractable for the $F_d^2$ class considered above.
Using again \eqref{eq:tay2} for all $x\in B(0,r_d), x\not=0$, we obtain

\begin{align}
\notag\int_{B(0,r_d)}|f(x)-\widetilde T_k f(x)|dx&\le \int_{B(0,r_d)}\frac{1}{k!}\|x\|_2^{k+1}\int_0^1(1-t)^k|(\partial_\nu^{k+1}f)(tx)|dtdx\\
\notag&= \frac{1}{k!}\int_{0<t<1} \int_{0< \|y\|_2< tr_d}\|y/t\|_2^{k+1}(1-t)^{k}|(\partial_\nu^{k+1}f)(y)|t^{-d}dydt\\
\notag&=\frac{1}{k!}\int_{B(0,r_d)}\|y\|_2^{k+1} \cdot |(\partial_\nu^{k+1}f)(y)|\cdot\int_{\|y\|_2/r_d}^1 t^{-(k+1)}(1-t)^k t^{-d}dtdy\\
\notag&\le \frac{1}{k!}\int_{B(0,r_d)}\|y\|_2^{k+1} \cdot |(\partial_\nu^{k+1}f)(y)|\cdot\int_{\|y\|_2/r_d}^1 t^{-(k+d+1)}dtdy\\
\label{eq:center}&\le\frac{1}{k!}\int_{B(0,r_d)}\|y\|_2^{k+1} \cdot |(\partial_\nu^{k+1}f)(y)| \cdot \frac{1}{k+d}\cdot \Bigl(\frac{\|y\|_2}{r_d}\Bigr)^{-k-d}dy\\
\notag&= \frac{1}{k!}\cdot \frac{r_d^{k+d}}{k+d}\int_{B(0,r_d)}\|y\|_2^{1-d}\cdot |(\partial^{k+1}_\nu f)(y)|dy\\
\notag&= \frac{1}{k!}\cdot \frac{r_d^{k+d}}{k+d}\int_{0}^{r_d}r^{1-d}\int_{r{\mathbb S}^{d-1}}|(\partial^{k+1}_\nu f)(y)|d\sigma(y)dr,
\end{align}
where $\sigma$ is the $d-1$ dimensional Hausdorff measure in $\R^d$ and ${\mathbb S}^{d-1}=\{x\in\R^d:\|x\|_2=1\}$ is the unit sphere in $\R^d.$
If we denote by $\omega_{d-1}$ the surface area of ${\mathbb S}^{d-1}$, i.e. $\omega_{d-1}=\sigma({\mathbb S}^{d-1})$, and by
$$
S(\partial_\nu^{k+1}f,r)=\frac{1}{\omega_{d-1}r^{d-1}}\int_{r {\mathbb S}^{d-1}}|\partial_\nu^{k+1}f(y)|d\sigma(y)
$$
the averages of $|\partial_\nu^{k+1}f|$ on the sphere $r{\mathbb S}^{d-1}$, we obtain

\begin{align*}
\int_{B(0,r_d)}|f(x)-\widetilde T_k f(x)|dx&\le \frac{1}{k!}\cdot \frac{r_d^{k+d}}{k+d}\cdot \omega_{d-1}\int_0^{r_d} S(\partial_\nu^{k+1}f,r) dr.
\end{align*}
Assuming finally, that the volume of the $B(0,r_d)$ is equal to $\omega_{d-1}r_d^{d}/d=1$, we may further reduce this to
\begin{align*}
\int_{B(0,r_d)}|f(x)-\widetilde T_k f(x)|dx&\le \frac{r_d^k}{k!} \int_0^{r_d} S(\partial_\nu^{k+1}f,r) dr,
\end{align*}
which is smaller than $r_d^k/k!$ for every $f\in F_d^3$, cf. \eqref{eq:Fd3},
or smaller than $r_d^{k+1}/k!$ for every $f\in F_d^4$, where
$$
F_d^4=\Bigl\{f\in C^\infty(B(0,r_d)):\sup_{k\in\N_0} \sup_{0<r\le r_d }S(\partial_\nu^{k+1}f,r)\le 1\Bigr\}.
$$

The proof of weak tractability follows in both cases from these estimates exactly as in Section \ref{sec:eucl}.

\section{Extensions and closing remarks}

The aim of this paper was to discuss the use of Taylor's theorem for tractability of approximation of infinitely differentiable functions.
Therefore, we have restricted ourselves in Theorem \ref{thm1.1} to the classes $F_d^1, F_d^2, F_d^3$ and simple domains as balls and cubes.
Nevertheless, the method used here may be directly generalized also to other situations.

{\bf Remark (Possible extensions):}
{\rm (i)} The point (i) of Theorem \ref{thm1.1}, which deals with uniform approximation on the cube for the class $F_d^1$, holds also for the classes
\begin{equation}\label{eq:Fd1c}
F_d^1(c)=\Bigl\{f\in C^\infty([-1/2,1/2]^d):\sum_{|\beta|=k} \frac{\|D^\beta f\|_\infty}{\beta!}\le c^k\quad\text{for all}\quad k\in\N_0\Bigr\},
\end{equation}
where $c<2$ is a fixed number. The analysis done in Section \ref{Sec2.1} applies literally also to this setting.\\
{\rm (ii)} The statements of points (ii) and (iii) of Theorem \ref{thm1.1} hold true also when $F_d^2$ and $F_d^3$ are replaced by
\begin{equation}\label{eq:Fd2c}
F_d^2(c)=\{f\in C^\infty(B(0,r_d)):\|\partial^{k}_\nu f\|_\infty\le c^k\quad\text{for all}\quad k\in\N_0\}
\end{equation}
or
\begin{equation}\label{eq:Fd3c}
F_d^3(c)=\Bigl\{f\in C^\infty(B(0,r_d)): \int_0^{r_d} S(\partial_\nu^{k}f,r) dr\le c^k\quad\text{for all}\quad k\in\N_0\Bigr\},
\end{equation}
respectively. Here, $c<\infty$ is arbitrary. Actually, one observes that even $c=d^{1/2-\delta}$ for some fixed $\delta>0$ is still admissible.
Again, also in these cases, the analysis done in Section \ref{sec:eucl} or Section \ref{sec:numer} applies.\\
{\rm (iii)} The analysis of Section \ref{sec:eucl} was done for a sequence of balls $B(0,r_d)\subset \R^d$. Although we concentrated in the very end on the
case, when the volume of $B(0,r_d)$ is equal to one (i.e. $r_d\approx\sqrt{d}$), the same calculation applies also to the case when $r_d\le C d^{1-\delta}$
for two universal constants $C>0$ and $\delta>0.$ Furthermore, the same is true for star-shaped subsets of such balls with zero in their center. Therefore, the statement applies
also to unit cube $[-1/2,1/2]^d$.\\
{\rm (iv)} When trying to generalize the analysis of Section \ref{sec:numer} to other domains, we encounter several problems. The crucial calculation \eqref{eq:center}
made a heavy use of spherical coordinates and they were also used in the definition of $S(\partial_{\nu}^{k+1}f,r)$. Although these obstacles can be overcome by
measuring the size of $\partial^{k+1}_\nu f$ in a certain weighted space, we avoid the technicalities and do not give the details.

{\bf Remark (Open problems):} (i) We have provided only the upper bounds on $n(\varepsilon,d)$, which in turn led to tractability results for the classes $F_d^1, F_d^2$, and $F_d^3$, respectively.
It would be interesting to know, if these results are optimal, i.e. to show that these results can not be improved. We leave this as an open problem.\\
(ii) We studied only algorithms using the function values of $f$. Of course, the presented tractability results also apply to the larger class of algorithms using
arbitrary linear functionals. Nevertheless, we leave as an open problem if such algorithms may achieve a better art of tractability than the one presented here.\\
(iii) One could also modify the classes under consideration in a way used recently in \cite{HNUW2}. This approach uses a sequence $L=(L_d)_{d\in \N}$ to define, for example,
\begin{equation*}
F_d^2(L)=\{f\in C^\infty(B(0,r_d)):\sup_{k\in \N_0} \|\partial^{k}_\nu f\|_\infty\le L_d\}.
\end{equation*}
Of course, the calculations given above could be to some extent transferred also to this setting, and we could really prove similar results for sequences which do not grow too quickly.
Unfortunately, due to the lack of lower estimates, we could not hope for being able to \emph{characterize} sequences $L$, for which weak (or quasi-polynomial) tractability still holds.

{\bf Acknowledgment:} I would like to thank to Aicke Hinrichs, Erich Novak, Mario Ullrich, Markus Weimar and Henryk Wo\'zniakowski for useful comments, remarks and hints
to an earlier version of the manuscript. Furthermore, I would like to thank to anonymous referees for their valuable comments, which helped to improve the presentation of the paper.

\end{document}